\documentclass[12pt, reqno]{amsart}
\usepackage{amssymb}
\usepackage{amsfonts}
\usepackage{amssymb,amsmath,amscd}
\usepackage[colorlinks=true]{hyperref}
\usepackage{mathrsfs}
\newtheorem{Theorem}{\indent Theorem}[section]

\newtheorem{Lemma}[Theorem]{\indent Lemma}

\theoremstyle{remark}
\newtheorem{Remark}{Remark}

\begin{document}
\centerline{\bf On an exponential sum  related to the M\"{o}bius function
}

\bigskip
\centerline{\small Wei Zhang}
\centerline{School of Mathematics and Statistics, Henan University}
\centerline{zhangweimath@126.com}
\centerline{https://orcid.org/0000-0002-5505-2671}
\bigskip

\textbf{Abstract}\
Let $\mu(n)$ be the M\"{o}bius function and $e(\alpha)=e^{2\pi i\alpha}$. In this paper, we study upper bounds of the classical sum
\[
S(x,\alpha):=\sum_{1\leq n\leq x}\mu(n)e(\alpha n).
\]
 We can improve some classical results of Baker and Harman \cite{BH}.

\medskip


\bigskip
\bigskip
\numberwithin{equation}{section}

\section{Introduction}
Let $S(x,\alpha)=\sum_{1\leq n\leq x}\mu(n)e(\alpha n).$   Davenport \cite{Da}  proved that
\[\max_{\alpha\in[0,1]}|S(x,\alpha)|\ll_{A} x(\log x)^{-A}
\]
for any $A>0.$  Hajela and  Smith  \cite{HS} gave the following conditional bounds for $S(x,\alpha)$. Suppose that for every Dirichlet character $\chi(\textup{mod}\ l)$, $L(s,\chi)$ has no zeros in the half-plane $\sigma>a$. Then \[\max_{\alpha\in[0,1]}|S(x,\alpha)|\ll x^{(a+2)/3+\varepsilon}\]
 for every $\varepsilon>0.$ Subsequently, Baker and Harman \cite{BH} sharpen the estimate to
 \[
 \max_{\alpha\in[0,1]}|S(x,\alpha)|\ll x^{b+\varepsilon},
 \]
 where
 \begin{align*}
 b=\begin{cases}
 a+1/4 &\textup{for}\ 1/2\leq a <11/20,\\
 4/5    &\textup{for}\ 11/20\leq a <3/5,\\
 (a+1)/2 &\textup{for}\ 3/5\leq a <1.
 \end{cases}
 \end{align*}
  In proving the above bound, in \cite{BH}, the authors   established that for any rational number $r/q$ with $(r,q)=1,$
\[
S(x,\alpha)\ll x^{a+\varepsilon}q^{1/2}(1+x|\alpha-r/q|)^{1/2}.
\]
The corresponding estimate made implicit by Hajela and Smith \cite{HS} is
\[
S(x,\alpha)\ll x^{a+\varepsilon}q^{1/2}(1+x|\alpha-r/q|).
\]
And the improvement comes by exploiting cancellations that occur in auxiliary sums and integrals.

  The aim of this paper is to give further improvement of  the above bound for certain range of $a$.

The idea is as follows: the
weak Generalized Riemann Hypothesis that all Dirichlet L-functions $L(s,\chi)$  have
no zeros in the half-plane $\sigma>a$  with $a\in[1/2,1)$
and the zero density results of Ingham \cite{In} imply  that
\begin{align}\label{In}
\sum_{\chi(\textup{mod}\ l)}N(\sigma,T,\chi)\ll (lT)^{(3-3\sigma)/(2-\sigma)+\varepsilon},
\end{align}
for $1/2\leq \sigma<a$ and
$0$ for $a\leq \sigma<1,$
where
\begin{align*}
N(\sigma,T,\chi):=\sum_{L(s,\chi)=0
\atop \sigma\leq\beta\leq1
,\ 10\leq|\gamma|\leq T}1.
\end{align*}
Hence, by using such a basic observation, we can obtain the following result. We can improve the classical result of Baker and Harman \cite{BH} for $1/2<a<4/7.$
\begin{Theorem}
Let $l$ be a positive integer and $\alpha$ be a real number.
Suppose that for every Dirichlet character $\chi(\textup{mod}\ l)$, $L(s,\chi)$ has no zeros in the half-plane $\sigma>a$. Then with the notions above, we have
 \[
 \max_{\alpha\in[0,1]}|S(x,\alpha)|\ll x^{b+\varepsilon},
 \]
 where
 $a\in[1/2,4/7]$ and
 \begin{align*}
 b=\frac{8a-7a^{2}}{4-2a}.
 \end{align*}
\end{Theorem}
\begin{Remark}
For example, for $1/2\leq a \leq 4/7,$ we have
\[
1/4+a\geq\frac{8a-7a^{2}}{4-2a}
\]
with equality happening if and only if $a=1/2,$ and
\[
4/5\geq\frac{8a-7a^{2}}{4-2a},
\]
with equality happening if and only if $a=4/7.$
Hence we have a much better result.
\end{Remark}

On the other hand, our idea also implies the following conditional result.
\begin{Theorem}
Let $q$ be a positive integer and $\alpha$ be a real number.
Suppose that for every Dirichlet character $\chi(\textup{mod}\ l)$, $L(s,\chi)$ has no zeros in the half-plane $\sigma>a$. Moreover, we also assume that
\begin{align*}
\sum_{\chi(\textup{mod}\ l)}N(\sigma,T,\chi)\ll (lT)^{2-2\sigma+\varepsilon},
\end{align*}
where
\begin{align*}
N(\sigma,T,\chi):=\sum_{L(s,\chi)=0
\atop \sigma\leq\beta\leq1
,\ 10\leq|\gamma|\leq T}1.
\end{align*}
Then with the notions above, we have
 \[
 \max_{\alpha\in[0,1]}|S(x,\alpha)|\ll x^{b+\varepsilon},
 \]
 where
 \begin{align*}
 b=\frac{1+a}{2}.
 \end{align*}
\end{Theorem}
\section{Proof of the main results}
Before the proof we will quote the following upper bound of $S(x,\alpha),$ which is related to the zero density results and has been used to deal with the cases of almost all in \cite{MS}.
\begin{Lemma}[See \cite{MS}]\label{MS}
Let $q$ be a positive integer and $\alpha$ be a real number. Then for any rational number $r/q$ with $(r,q)=1$ and for any fixed $\theta$ satisfying $1/2+\delta\leq \theta\leq 1-\delta$ (with $\delta$ being any small positive constant), we have
\begin{align*}
\begin{split}
S(x,\alpha) &\ll q^{-1/2}x^{1+\varepsilon}+
\sum_{d|q}(q/d)^{1/2}(\phi(q/d))^{-1}(1+x|\alpha-r/q|)\\
& \times
\left(
\max_{1\leq T\leq (x/d)^{3}}T^{-1}(qT)^{\varepsilon}
\int_{1/2}^{1}
\left(
T\phi(q/d)(x/d)^{\theta+\varepsilon}+
(x/d)^{\sigma+\varepsilon}\sum_{\chi(\textup{mod}\ q/d)}N(\sigma,T,\chi\chi_{d})
\right)d\sigma
\right),
\end{split}
\end{align*}
where $\phi(n)$ is the Euler function and
\begin{align*}
N(\sigma,T,\chi):=\sum_{L(s,\chi)=0
\atop \sigma\leq\beta\leq1
,\ 10\leq|\gamma|\leq T}1.
\end{align*}
\end{Lemma}
On the other hand, we also need the following result of Baker and Harman \cite{BH}.
\begin{Lemma}[See \cite{BH}]\label{BH}
Suppose that for every Dirichlet character $\chi(\textup{mod}\ l)$, $L(s,\chi)$ has no zeros in the half-plane $\sigma>a$. Then with the notions above, for any rational number $r/q$ with $(r,q)=1,$
\[
S(x,\alpha)\ll x^{a+\varepsilon}q^{1/2}(1+x|\alpha-r/q|)^{1/2}.
\]
\end{Lemma}

We now apply Dirichlet's theorem to obtain, for any real number $\alpha$, a rational
number $r/q$ with $(r,q)=1,$ $1\leq q \leq x^{a}$ such that
\[
\left|\alpha-\frac{r}{q}\right|<\frac{1}{qx^{a}}.
\]
By Lemma \ref{BH}, for $q\leq x^{1-a},$ we have
\[
S(x,\alpha)\ll x^{a+\varepsilon}q^{1/2}(1+x|\alpha-r/q|)^{1/2}
\ll x^{(1+a)/2+\varepsilon}.
\]
Then by Lemma \ref{MS} with $\theta=1/2+\varepsilon$, the hypothesis such that  for every Dirichlet character $\chi(\textup{mod} \ l)$, $L(s,\chi)$ has no zeros in the half-plane $\sigma>a$, we have
\begin{align*}
S(x,\alpha)&\ll q^{-1/2}x^{1+\varepsilon}+
x^{1/2+\varepsilon}q^{1/2}(1+x|\alpha-r/q|)
\\&+
q^{-1/2}x^{a+\varepsilon}q^{3(1-a)/(2-a)}(1+x|\alpha-r/q|).
\end{align*}

It is worth pointing out that by the hypothesis,    in Lemma  \ref{MS} , the range of integration for the second term inside the maximum
over $T$ is reduced to $[1/2, a]$ (for the possible non-zeros in $N(\sigma,T,\chi\chi_{d})$, Ingham's
zero density result (\ref{In}) is applied).

Then for $x^{1-a}\leq q\leq x^{a},$ we have
\[
x|\alpha-r/q|\leq 1
\]
and
\[
x^{1/2+\varepsilon}q^{1/2}+x^{1+\varepsilon}q^{-1/2}\ll x^{(1+a)/2+\varepsilon}.
\]
Hence,  for $a\in[1/2,3/4],$ we  can obtain
\begin{align*}
S(x,\alpha)&\ll x^{1/2+\varepsilon}q^{1/2}+x^{1+\varepsilon}q^{-1/2}
+q^{-1/2}x^{a+\varepsilon}q^{3(1-a)/(2-a)}\\
& \ll x^{(8a-7a^{2})/(4-2a)+\varepsilon}.
\end{align*}
This completes the proof.

\bigskip
 {\bf Acknowledgements} The author would like to thank the the referee who gives some  detailed corrections and suggestions.

\address{Wei Zhang\\ School of Mathematics and Statistics\\
               Henan University\\
               Kaifeng  475004, Henan\\
               China}
\email{zhangweimath@126.com}

\end{document}